\documentclass[]{amsart}

\usepackage{latexsym}
\usepackage{amssymb,amsmath,amsopn}
\usepackage[dvips]{graphicx}   
\usepackage{wrapfig}
\usepackage{subfig}
\usepackage{color,epsfig}      
\usepackage{url}
\usepackage{amsbsy}

\newtheorem{statement}{Statement}[]

\newtheorem{defi}{Definition}[]

\theoremstyle{definition}

\newtheorem{conjecture}{Conjecture}[]

\begin{document}
\title[On the constructibility of the axes of an ellipsoid]{On the constructibility of the axes of an ellipsoid.\\[1ex]
  \footnotesize\mdseries The construction of Chasles in practice}
\author[\'A. G.Horv\'ath]{\'Akos G.Horv\'ath}
\address {Department of Geometry \\
Budapest University of Technology and Economics\\
H-1521 Budapest\\
Hungary}
\email{ghorvath@math.bme.hu}
\author[I. Prok]{Istv\'an Prok}
\address {Department of Geometry \\
Budapest University of Technology and Economics\\
H-1521 Budapest\\
Hungary}
\email{prok@math.bme.hu}

\dedicatory{Dedicated to the memory of the teaching of Descriptive Geometry in BME Faculty of Mechanical Engineering}

\date{April, 2017}

\begin{abstract}
In this paper we discuss Chasles's construction on ellipsoid to draw the semi-axes from a complete system of conjugate diameters. We prove that there is such situation when the construction is not planar (the needed points cannot be constructed with compasses and ruler) and give some others in which the construction is planar.
\end{abstract}

\maketitle

\section{Introduction}
Who interested in conics know the construction of Rytz, namely a construction which determines the axes of an ellipse from a pair of conjugate diameters. Contrary to this very few mathematicians know that in the first half of the nineteens century Chasles (see in \cite{chasles}) gave a construction in space to the analogous problem. The only reference which we found in English is also in an old book written by Salmon (see in \cite{salmon}) on the analytic geometry of the three-dimensional space. In \cite{gho:conf} the author extracted this analytic method to prove some results on the quadric of $n$-dimensional space. However, in practice this construction cannot be done using compasses and ruler only as we will see in present paper (Statement \ref{st:nonconstr}). On the other hand those steps of the construction which are planar constructions we can draw by descriptive geometry (see the figures in this article).

\section{Basic properties}

The proof of the statements of this section can be found in \cite{gho:conf}. In this paper we use a Cartesian coordinate system in the three space (or in the $n$-space). A non-degenerated central surface of second order is called by \emph{quadric} if it has the following canonical form:
\begin{equation}\label{eq:quadric}
\sum\limits_{i=1}^n\varepsilon_i\frac{x_i^2}{(a_i)^2}=1,
\end{equation}
with $0< a_n\leq a_{n-1}\leq \ldots \leq a_1$ and $\varepsilon_i\in\{\pm 1\}$. Clearly, if $x'$ is any point of the quadric above a normal vector at $x'$ is equal to
\begin{equation}\label{eq:normal}
n(x')=\left(\varepsilon_1\frac{x'_1}{(a_1)^2},\ldots, \varepsilon_n\frac{x'_n}{(a_n)^2}\right)^T \mbox{ with norm square } |n(x')|^2=\sum\limits_{i=1}^{n}\frac{(x'_i)^2}{(a_i)^4}
\end{equation}
and its tangent hyperplane at this point is the set of points $x$ with
$$
0=\langle n(x'),x-x'\rangle=\sum\limits_{i=1}^{n}\varepsilon_i\frac{x_ix'_i-{x'_i}^2}{(a_i)^2}=\sum\limits_{i=1}^{n}\varepsilon_i\frac{x_ix'_i}{(a_i)^2}-1
$$
We use the name \emph{ellipsoid} for a quadric in that case when for all $i$ hold $\varepsilon_i=1$. If for all $i$ the number $\varepsilon$ is $-1$, there is no real point of the quadric.

\begin{defi}\label{def:conjdiam} Two diameters of the ellipsoid are \emph{conjugate}, if they are conjugate diameters in the $2$-dimensional ellipse is the intersection of the quadric by the $2$-plane through the diameters.
\end{defi}

The line parallel to the first diameter and through the end point of the second one is a tangent line of the quadric at this point and vica versa. Analytically we have a condition for the conjugate directions $e=(e_1,\ldots, e_n)^T$ and $f=(f_1,\ldots,f_n)^T$. Let $p_e, p_f$ be such real numbers that $p_ee,p_ff$ are the points of the quadric. Then the assumption above says that
\begin{equation}\label{eq:conjugatediam}
0=\frac{1}{p_e}\langle n(p_ee),f\rangle=\sum\limits_{i=1}^{n}\varepsilon_i\frac{e_if_i}{(a_i)^2}
\end{equation}
A complete system of conjugate semi-diameters $\{x^1,\ldots,x^n\}$ has $n$ pairwise conjugate elements. For an ellipsoid the above conditions can be given in matrix form, too. Let denote by $X$, $X_A$ and $A$ the respective matrices
$$
X=\left(
  \begin{array}{ccc}
    x^1_1 & \cdots & x^n_1 \\
    \vdots & \cdots & \vdots \\
    x^1_n & \cdots & x^n_n \\
  \end{array}
\right), \quad
X_A=\left(
      \begin{array}{ccc}
        \frac{x^1_1}{a_1} & \cdots & \frac{x^1_n}{a_n} \\
        \vdots & \cdots & \vdots \\
        \frac{x^n_1}{a_1} & \cdots & \frac{x^n_n}{a_n} \\
      \end{array}
    \right)
\mbox{ and }
A=\left(
    \begin{array}{ccccc}
      a_1 & 0 & \cdots & 0 & 0\\
      0 & a_2 & 0 &\cdots & 0 \\
      0 & 0& a_3 & \cdots & 0 \\
      \vdots & \vdots &\vdots & \vdots &\vdots \\
      0 & 0 & \cdots & 0 & a_n \\
    \end{array}
  \right).
$$
Clearly $X_A\cdot A=X$ and $X_A$ by the above assumption is an orthogonal matrix. From this follows that
$$
X^TX=\left(X_A\cdot A\right)^T\left(X_A\cdot A\right)=A^TA=A^2
$$
implying that the respective traces of the two sides are equal to each other. This means that
\begin{equation}\label{eq:squaresumofthediam}
|x^1|^2+|x^2|^2+\ldots +|x^n|^2=(a_1)^2+\ldots +(a_n)^2.
\end{equation}
\emph{The sum of squares of a complete system of conjugate semi-diameters is a constant}. The orthogonality of the matrix $X_A$ immediately implies that its determinant is $1$, hence we also have:
\begin{equation}\label{eq:volumeofpar}
\det X =\det (X_A\cdot A)=\det A=a_1\cdots a_n,
\end{equation}
meaning that \emph{the volume of the parallelepiped spanned by the semi-diameters of a complete conjugate system is a constant}.

\begin{defi}\label{def:conjpoint}
Two points $x$ and $y$ are \emph{conjugate with respect to the ellipsoid} if holds the equality $x^TA^{-2}y=1$ with the diagonal matrix $A^{-2}$ containing the respective reciprocal values of the square of the semi-axes in its diagonal.
\end{defi}
Let $\langle H, x\rangle=\sum\limits_{i=1}^nh_ix_i=1$ is a given hyperplane not containing the origin. The \emph{pole} of this hyperplane with respect to the ellipsoid is such a point of the space which is conjugate to an arbitrary point of the hyperplane. If $\xi$ is its pole with respect to the ellipsoid then for a point $x\in H$ we have two equations for all $x\in H$:
$$
\langle H, x\rangle=1 \mbox{ and } 1=x^TA^{-2}\xi=\sum\limits_{i=1}^n \frac{x_i\xi_i}{(a_i)^2}.
$$
Now immediately follows that
\begin{equation}\label{eq:poleofaplane}
\xi_i=h_i(a_i)^2 \quad i=1,\ldots, n.
\end{equation}

Another important class of surfaces of second order the class of \emph{cones} of second order defined by the elements of a pencil of lines envelop a conic. If its apex is the origin it can be written in the (canonical) form
\begin{equation}\label{eq:cone}
\sum\limits_{i=1}^n\varepsilon_i\frac{x_i^2}{(a_i)^2}=0,
\end{equation}
with $a_n\leq a_{n-1}\leq \ldots \leq a_1$ and $\varepsilon_i\in\{\pm 1\}$. The \emph{semi-axes} of this cone are the pairwise orthogonal segments on the corresponding coordinate axes with respective lengths $a_1,\ldots,a_n$. Clearly the origin is the only point of a cone if the signs $\varepsilon_i$ are equal to each other.  In the other cases, the cones contain lines through the origin. The intersection of a proper cone with an affine hyperplane (not going through the origin) is a quadric of dimension $n-1$.

\begin{defi}\label{def:confcon}
Consider an ellipsoid in canonical form with squared semi-diameters $(a_i)^2$ $i=1,\ldots, n$. The equalities
\begin{equation}\label{eq:confocals}
\sum\limits_{i=1}^n\frac{x_i^2}{(a_i)^2-\lambda}=1,
\end{equation}
with different values of $\lambda$ define those quadrics which are \emph{confocal to the given ellipsoid}. We assume in this paper that $0<a_n<a_{n-1}<\ldots <a_1$ and $\lambda\in\mathbb{R}$ moreover we denote them by $\mathcal{C}(\lambda)$.
\end{defi}

First we collect some important properties of the confocal quadrics in general.

\begin{enumerate}
\item[{\bf I:}] \emph{Any point $x'=(x'_1,\ldots, x'_n)^T$ with non-zero coordinates of the space belongs to a system of confocal (to the given one) quadrics of cardinality $n$}. In fact, we can arrange the original equation to the form
\begin{equation}\label{eq:flambdafunc}
0=\sum\limits_{i=1}^n{x'_i}^2\prod\limits_{j\ne i}((a_j)^2-\lambda)-\prod\limits_{i=1}^n((a_i)^2-\lambda)=:f(\lambda).
\end{equation}
Since $f(-\infty)=-\infty ,\, f((a_n)^2)>0, \, f((a_{n-1})^2)<0, \ldots $ by Role's theorem there are $n$ distinct roots $\lambda^j$ for $\lambda$. For every solution we have a system of semi-axes denoted by $\{a_i^j \, i=1,\ldots n\}$ where $j=1,\ldots ,n$ (if we choose a point on the ellipsoid then a solution determines the original system of axes $a_i$, we can assume that $a_i^1:=a_i$). Note that even the value $a_i^j$ is the length of that semi-axis of the $j$-th system which corresponds to the $i$-th coordinate axis, with respect to the fixed distance of the lengths of the corresponding semi-axes in a confocal system these are again create a monotone decreasing sequence for a fixed $j$. On the other hand we also assume that the values $a_1^j$ gives a monotone decreasing sequence in $j$, namely $a_1^1 > a_1^2> \ldots > a_1^n$ holds. (We note that the greater semi axis of the given system of confocals is always
the greater semi axis of the only ellipsoid of the system. Hence there is no contradiction between our assumptions on $a_1^1$.)
\item[{\bf II:}] \emph{By the above squared semi-axes $\left(a_i^j\right)^2$ we can determine the coordinates of the common point $x'$}. It can be proved that
\begin{equation}\label{eq:coordwithconfoc}
{x'_i}^2=\frac{\prod\limits_{j=1}^n\left(a_i^j\right)^2} {\prod\limits_{\substack{j=1 \\ i\ne j}}^n\left(\left(a_i^i\right)^2-\left(a_j^i\right)^2\right)} \mbox{ for all } i.
\end{equation}

\item[{\bf III:}] \emph{The norm square of the point by the squared semi-axes is}
\begin{equation}\label{eq:lengthbytheconfaxes}
|x'|^2=\sum\limits_{i=1}^n{x'_i}^2=\sum\limits_{j=1}^n\left(a_j^j\right)^2,
\end{equation}
because of the equalities $\left(a_1^j\right)^2-\left((a_1)^2-(a_j)^2\right)=\left(a_1^j\right)^2-\left(\left(a_1^j\right)^2-\left(a_j^j\right)^2\right)=\left(a_j^j\right)^2$.

\item[{\bf IV:}] \emph{The confocal quadrics through the given point $x'$ are pairwise orthogonal at this point}. The equality for $j\ne k$
\begin{equation}\label{eq:onorthog}
0=\left(\left(a_1^k\right)^2-\left(a_1^j\right)^2\right)\sum\limits_{i=1}^{n}\frac{\left(x'_i\right)^2}{\left(a_i^j\right)^2\left(a_i^k\right)^2}
\end{equation}
is equivalent to the above statement. Therefore, where the $n$ confocal quadrics intersect, each tangent hyperplane cuts the other perpendicularly, and the tangent hyperplane to any one contains the normals to the other $n-1$.

\item[{\bf V:}] \emph{On the central sections of a quadric} it can be proved that if a plane be drawn through the centre parallel to any tangent hyperplane to a quadric, the axes of the section made by that plane are parallel to the normals to the $n-1$ confocals through the point of contact. Since the latter directions are pairwise orthogonal by {\bf IV}, it have to prove only that these directions are also conjugate with respect to the intersection quadric of dimension $n-1$.

\item[{\bf VI:}] \emph{The lengths of the semi-axes of the central section parallel to the tangent hyperplane at the point $x'$} are the respective lengths $\rho^1_j$ of the radius vector of the given quadric (which is the first confocal ones) parallel to the normal of the $j$-th confocal. From this
we get that for all $j\geq 2$ the squared length of the $j$-th semi-axis is equal to:
\begin{equation}\label{eq:radius}
\left(\rho^1_j\right)^2=\left(a_1^1\right)^2-\left(a_1^j\right)^2.
\end{equation}
\end{enumerate}
We can define a \emph{duality property} of confocal systems.
Denote by $(p^j)^2$ the reciprocal of the norm square of the $j$-th normal vector $n_j(x')$ at $x'$. Geometrically $p^j$ is the length of the foot-point of the perpendicular from the origin to the tangent hyperplane of the $j^{th}$ quadric. Now we have that for all $j$
$$
\left(p^j\right)^2=\frac{\prod\limits_{i=1}^n\left(a_i^j\right)^2}{\prod\limits_{k\ne j}\left|\left(a_1^j\right)^2-\left(a_1^k\right)^2\right|}.
$$
By the opportunity that we involve the sign of the squared value in the notation (as at (\ref{eq:coordwithconfoc})), we can write that for all $j$ we have
\begin{equation}\label{eq:onpj}
\left(p^j\right)^2=\frac{\prod\limits_{i=1}^n\left(a_i^j\right)^2}{\prod\limits_{k\ne j}\left(\left(a_1^j\right)^2-\left(a_1^k\right)^2\right)}
\end{equation}
defining a formula analogous that of (\ref{eq:coordwithconfoc}). It can be observed the symmetry which exists between these values for $p^i$, and the values already found in (\ref{eq:coordwithconfoc}) for $x'_i$. If the $n$ tangent hyperplanes had been taken as coordinate hyperplanes, $p^j$ would be the coordinates (with suitable signs) of the centre of the surface. So we have the following duality theorem:
\begin{statement}\label{st:duality}
With the point $x'$ as the centre $n$ confocals may be described having the $n$ tangent hyperplanes for principal planes and intersecting in the centre of the original system of surfaces. The axes of the new system of confocals are $a^1_1,a^2_1,\ldots,a^n_1$, $a^1_2,a^2_2,\ldots a^n_2$ $\ldots $ $a^1_n,a^2_n,\ldots a^n_n$. The $n$ tangent hyperplanes of the new (dual) system are the $n$ principal planes of the original system.
\end{statement}

Clearly, if $\lambda=a_k^2$ then we should consider only such points as a point of the quadric, which holds the assumption $x_k=0$ and the other coordinates satisfy the equation
\begin{equation}\label{eq:focalquadric}
\sum\limits_{\substack{i=1 \\ i\ne k}}^n\frac{x_i^2}{a_i^2-a_k^2}=1.
\end{equation}
Hence the above $n-1$-dimensional quadric can be considered as the limit of the pencil of confocal quadrics by the assumption $\lambda \rightarrow a_k$. Denote by $C_k$ this surface which we call the \emph{$k$-th focal quadric} of the pencil of quadric $\mathcal{C}(\lambda)$. First of all there are precisely $n-1$ focal surface with real points and one of them (when $\lambda=a_1^2$) does not contain real points.
It is clear that the $n$-th focal quadric $C_n$ is an ellipsoid.

The following statements on the focal quadrics can be used well.

\begin{statement}\label{st:axesoffocalcones}
Let $x'$ be any point of the space and denote by $C_k(x')$ the cones with apex $x'$ and generators through the focal quadric $C_k$. Then there are an orthogonal system of lines $l_1,\ldots, l_n$ through $x'$ which elements are common axes of the focal cones. The magnitudes of the axes corresponding to the same line are dependent from $k$, more precisely the signed and squared lengths of the semi-axes of $C_k(x')$ are $\{\left(a^i_k\right)^2 \,:\, i=1,\ldots,n\}$.
\end{statement}

\begin{statement}\label{st:lengthoftheintercepts}
$(n-1)$ cones having a common vertex $x$ envelope the $(n-1)$ focal quadrics (of distinct types). The length of the intercept made on one of their common edges by a hyperplane through the origin parallel to the tangent hyperplane to a confocal through $x$ is equal to the major semi-axis of the given confocal.
\end{statement}

\section{Construction the axes from a complete system of conjugate diameters}

\subsection{Once more again on the case of $n=2$.}\label{ssec:modrytz}

\begin{figure}[ht]
  \centering
  \includegraphics[width=15cm]{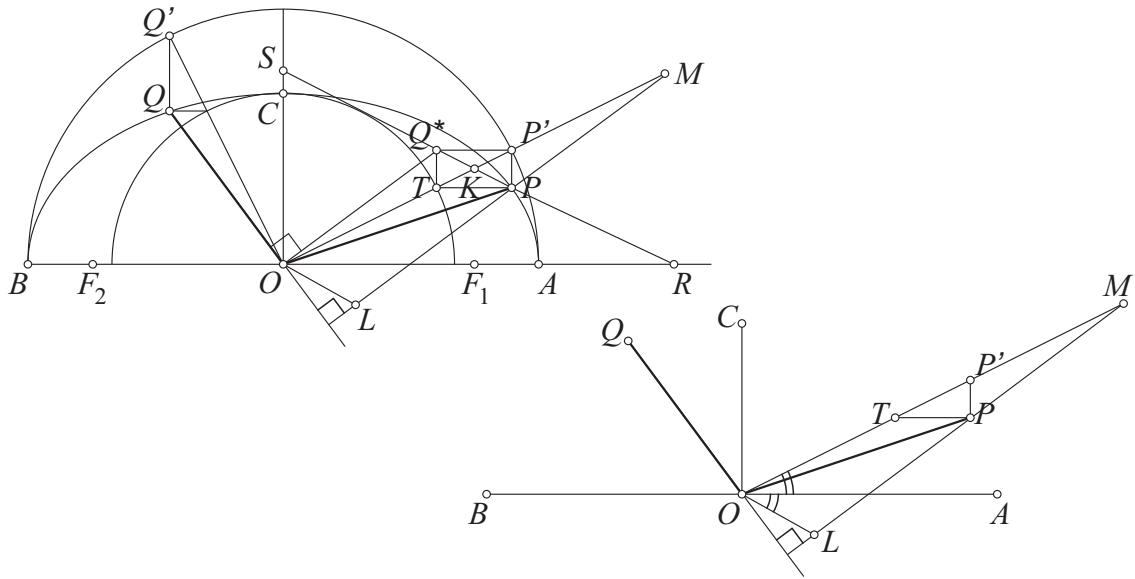}\\
  \caption{Conjugate diameters and axes }\label{fig:rytz}
\end{figure}
\begin{figure}[hb]
  \centering
  \includegraphics[height=9cm]{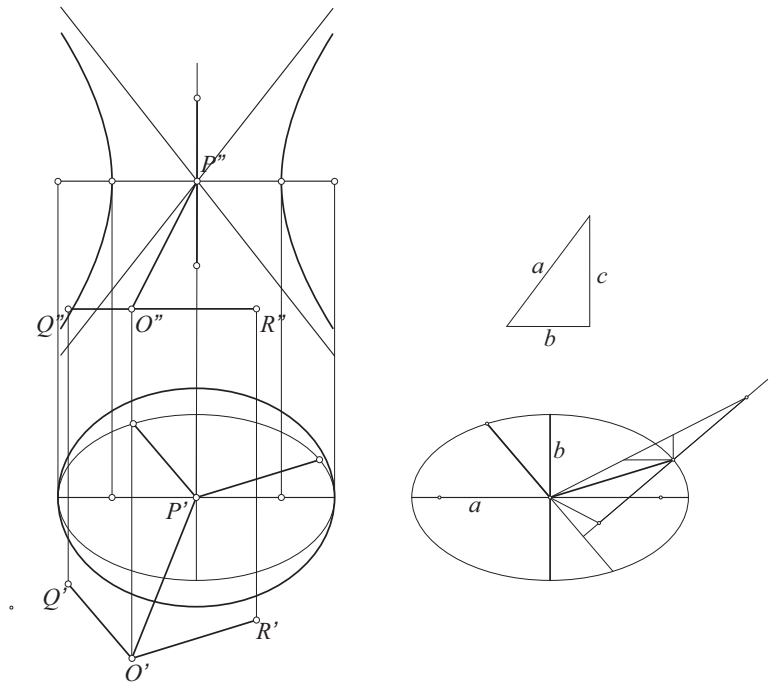}\\
  \caption{Construction of the focal conics}\label{fig:focalconics}
\end{figure}
The standard proof of the construction of Rytz is based on the so-called "two-circle figure" (see in Fig. \ref{fig:rytz}) in which we draw the incircle and the circumcircle of the ellipse. The same figure with a little modification enable to get another construction to solve this problem.
In fact, let the line $PM$ be the normal of the ellipse at $P$. This is perpendicular to $OQ$ hence it is parallel to $OQ^\star$. If $M$ is the intersection of this normal with the line $OK$ (as in the standard proof $Q^\star$ is the rotated copy of $OQ$ by $\pi/2$ and $K$ is the middle point of the segment $PQ^\star $), then the triangles $OTQ^\star$ and $MP'P$ are congruent to each other. From this we get that $PM$ is congruent to $OQ^\star$. Let $L$ be the reflected image of $M$ at $P$. Since $PL$ is parallel and equal to $OQ^\star$, $OL$ is is parallel and congruent to $Q^\star P$. Hence the axis $AB$ of the ellipse is the bisector of the angle $LOM\angle$. The other axis of the ellipse is perpendicular to $AB$ at $O$, and drawing parallels from $P$ to these axes we get the intersection points $T$ and $P'$ with the line $OM$, respectively. Clearly, the lengths of $OT$ and $OP'$ are equals to the lengths of the semi-axes, respectively. The construction has the following steps:

\begin{itemize}
\item Draw a perpendicular to the line $OQ$ from $P$ and determine the points $M$, $L$ on this line with property $|PM|=|PL|=|OQ|$, respectively.
\item Draw the bisectors of the angle $LOM\angle$, these are the axes of the ellipse.
\item Draw parallels to these bisectors from the point $P$ and determine the intersections of these lines with the line $OM$ ($T$, $P'$). The lengths of the semi-axes are the lengths of the segments $OT$, $OP'$, respectively.
\end{itemize}

Observe that $L$ and $M$ determine a pencil of confocal conics in the plane. From these pencils there are two conics an ellipse and a hyperbole through the point $O$. By the bisector property, the tangent lines of these conics at $O$ are the axes of the searching ellipse. If we visualize the union of lines $OM$ and $OL$ as a cone with apex $O$ and through the focal figure of the confocal system of conics determined by the points $M$ and $L$ then the axes of the ellipse are equals the axes of this cone. Hence this construction contains such concepts which analogously exist in every higher dimensions, too.

\subsection{The case of the three space.}

\begin{figure}[hb]
  \centering
  \includegraphics[height=12cm]{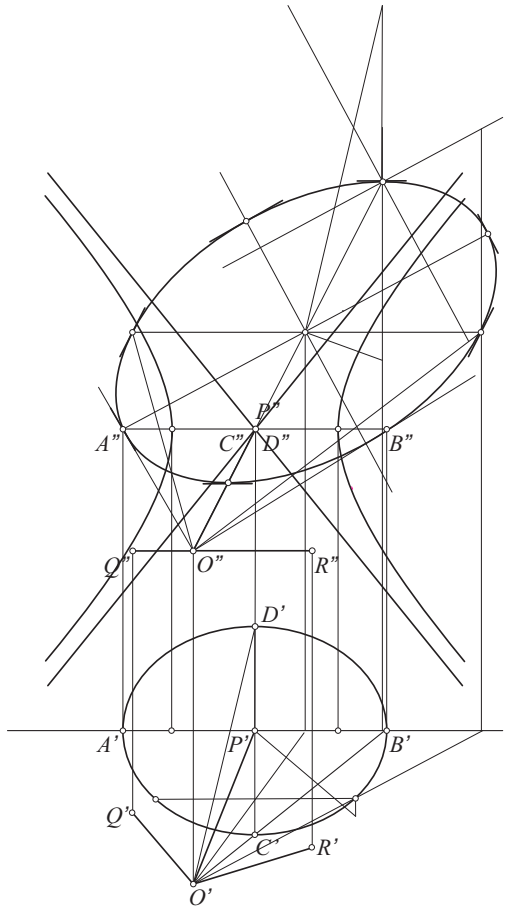}\\
  \caption{The image of the focal ellipse on the plane of the focal hyperbole}\label{fig:imoffocellipse}
\end{figure}

M. Chasles (see p.157 (182.$\S$) in \cite{salmon}, \cite{chasles}) observed that the above construction can be adopted to the three space, too.\footnote{In a more rigorous investigation shows that the "Chasles's construction" cannot be drawn by compasses and ruler since the construction of the common lines of two quadratic cones can be constructed if and only if the common point of two given conics can be constructed. This latter exercise can be solved some singular cases only, e.g. either if the two conics have a common focus or if they have common origin.} Let us see the theoretical "construction" of Chasles.

Consider the three pairwise conjugate semi-diameters $OP$,$OQ$ and $OR$, respectively.
\begin{itemize}
\item Determine the three mutually perpendicular line through $P$, one of its the normal of the ellipsoid (perpendicular to the plane $(OQR)$ at $P$) and the other two are parallel to the axes of the ellipse $\mathcal{E}$ which is the intersection of the plane $(OQR)$ by the given ellipsoid. (By the construction of Subsection \ref{ssec:modrytz} we can get these lines.)

\item They contains the axes of a confocal system of quadrics dual to the confocal system of the given ellipsoid. The duality statement  (Statement \ref{st:duality}) says that the axes of this systems can be given by the system $\{a_1^1,a_1^2,a_1^3\}$ as respective lengths of the semi-axes. The square of the semi-axes of the focal conics of this system are $\{(a_1^1)^2-(a_1^3)^2,(a_1^2)^2-(a_1^3)^2\}$; $\{(a_1^1)^2-(a_1^2)^2,(a_1^3)^2-(a_1^2)^2\}$ defining an ellipse and a hyperbole, respectively. The square of the axes of $\mathcal{E}$ are $(\rho_2^1)^2=(a_1^1)^2-(a_1^2)^2$ and $(\rho_3^1)^2=(a_1^1)^2-(a_1^3)^2$. These values are known from the construction of Subsection \ref{ssec:modrytz}, hence we can get also the squares of the lengths of the semi-axes of the focal conics.

\item The intersection of the two focal cones with apex $O$ is the union of the four common edges. The six planes determining by these edges intersect to each other in three mutual perpendicular lines which are the three common axes of these two confocal cones.\footnote{To prove this statement observe that the two focal cones are confocal quadrics and so their principal axes have common directions. The three principal axis intersects a general plane of intersection in three points which are form an autopolar triangle of this plane with respect to the two corresponding conics of intersection. The four common points of the two conic of intersections determines that quadrangle which diagonal points form the only autopolar triangle with respect to both of the conics.}  By Statement \ref{st:axesoffocalcones} these are the searched axes of the given ellipsoid, while the planes through $P$ and parallel to the principal planes of the ellipsoid cut off on these four lines parts equal in length to the semi-axes by Statement \ref{st:lengthoftheintercepts}.
\end{itemize}

The construction of the first two steps can be seen (after an appropriate transformation of the given diameters) on Fig. \ref{fig:focalconics}

\subsubsection{On the common lines of the cones which envelope the focal conics}

\begin{figure}[!ht]
     \subfloat[The lines of the axes\label{fig:axesinspace}]{%
       \includegraphics[width=0.45\textwidth]{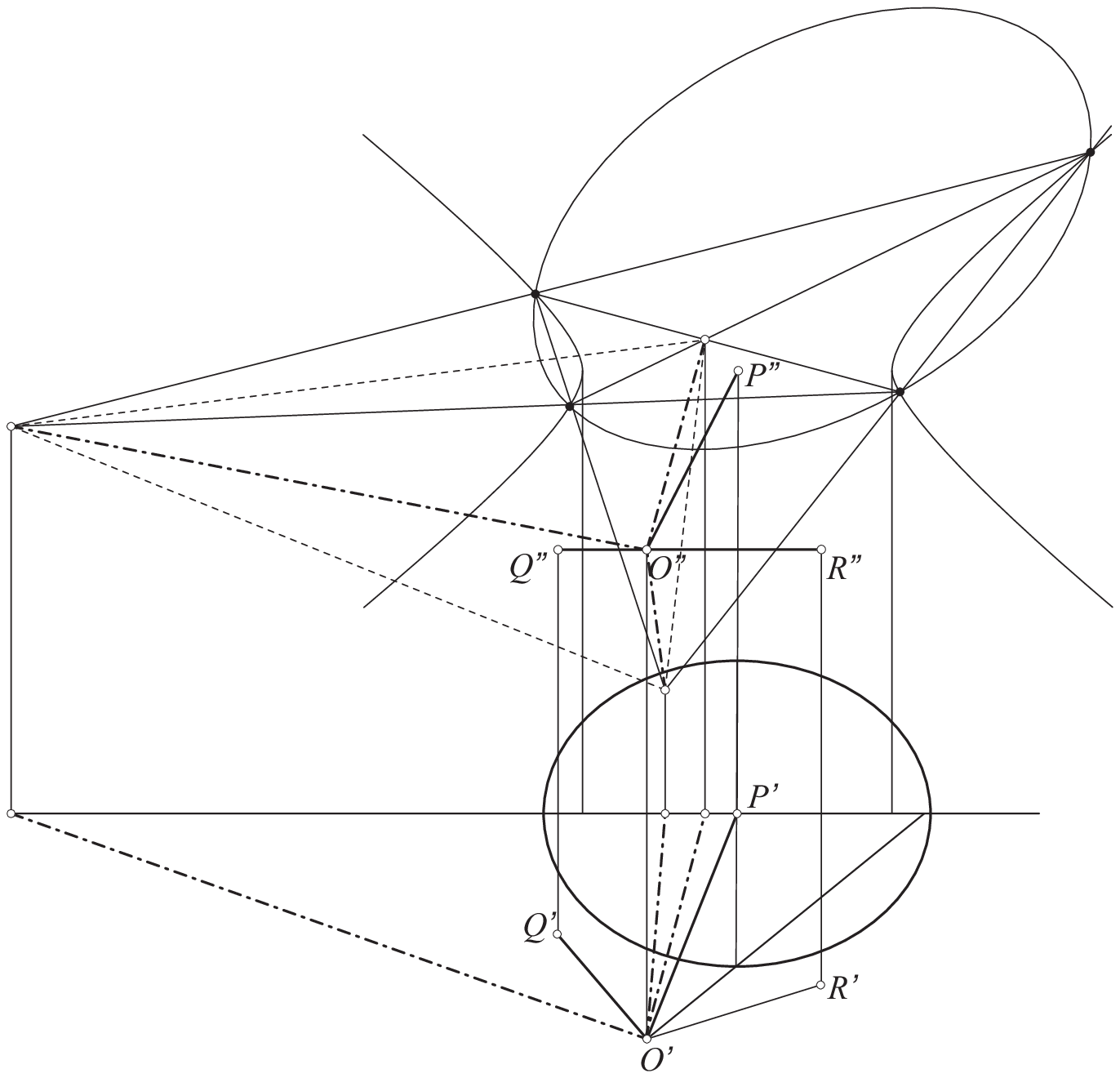}
     }
     \hfill
     \subfloat[Construction of the length of the major semi-axis $OA$\label{fig:lengthofaxes}]{%
       \includegraphics[width=0.45\textwidth]{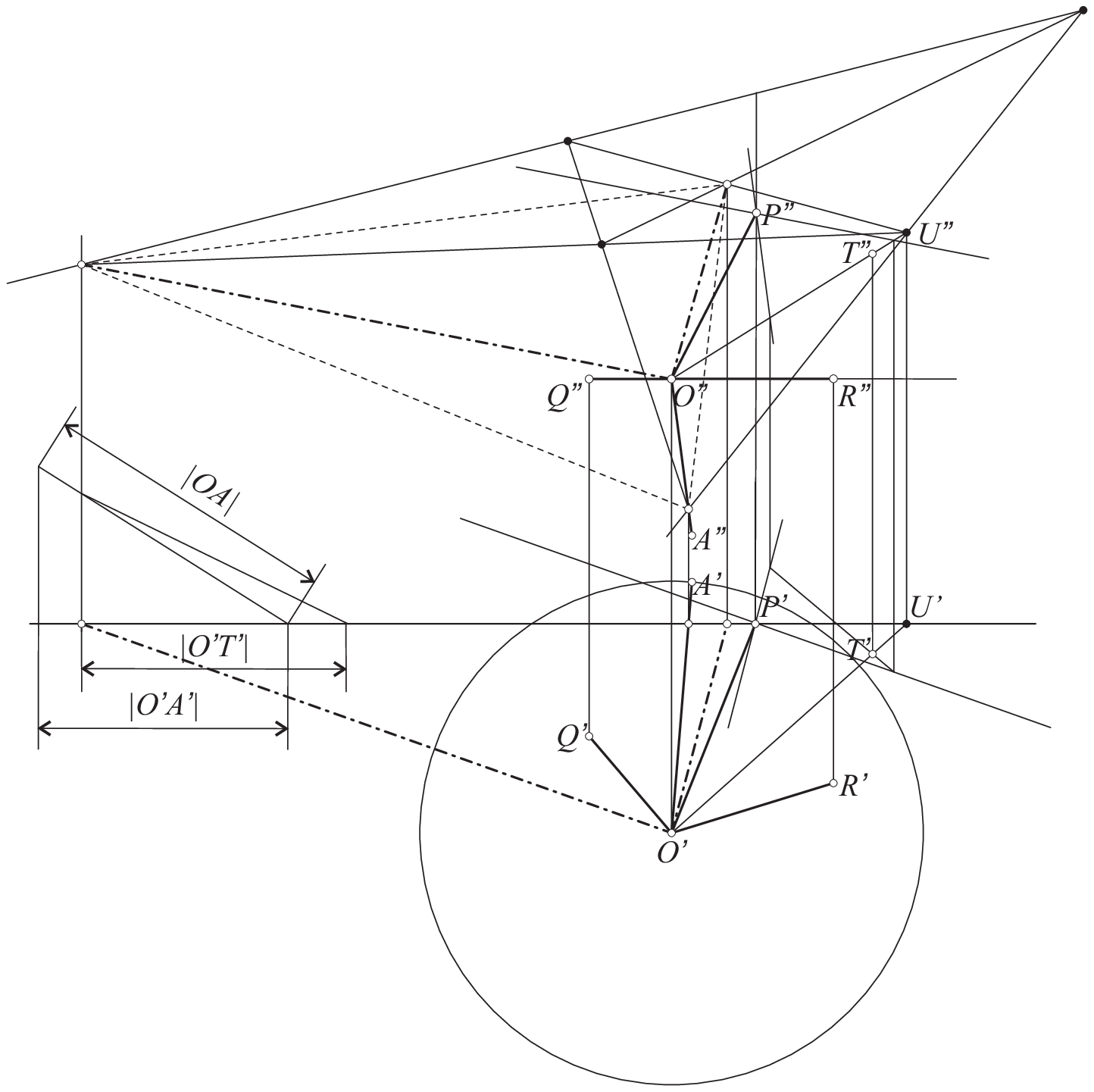}
     }
     \caption{The last step}
     \label{fig:thelaststep}
   \end{figure}

It is clear that in the last step of Chasles's construction we have to determine the common edges of the two focal cones. If the point $x'=O$ is given the common edges of the two cones intersects the plane of the focal hyperbole in the common points of the hyperbole and that conic $\mathcal{E}'$ which can be get from the focal ellipse by the central projection from $O$. The major axis $AB$ of the focal ellipse is fixed by this projection because it is on the intersection $m$ of the image and domain planes. We can construct the images of the tangent lines of the ellipse in its endpoint $A$ and $B$, and also the images $\overline{C}$ and $\overline{D}$ of the endpoints of the minor axis $CD$. Since the tangent lines in $C$ and $D$ are parallel to the line $m$ this is also true for their images, respectively. Hence $\overline{C}\overline{D}$ is a diameter of the image conic. Since the tangent lines at $\overline{C}$ and $\overline{D}$ are parallel to $m$, the midpoint of $\overline{C}\overline{D}$ is the midpoint of the diameter $\overline{E}\overline{F}$ of the image conic conjugate to $\overline{C}\overline{D}$, its direction is also parallel to $m$. Now we can construct that secant $EF$ of the focal ellipse which endpoints mapped to the endpoints of the image conic. The axes of $\mathcal{E}'$ (e.g. in the case when $\mathcal{E}'$ is an ellipse) can be get by the modified Rytz construction (see Fig. \ref{fig:imoffocellipse}). Finally, we would like to construct the common points of the focal hyperbole and the conic $\mathcal{E}'$. This step in general cannot be constructed using only straightedge and compasses. In the next section we discuss this problem a little bit more.
But if we have the four common edges of the focal cones we can construct the direction of the three axes and also the lengths of them. In Fig. \ref{fig:axesinspace}, Fig. \ref{fig:lengthofaxes} we finished the construction, assuming that we know the common points of the above conics, respectively.

\subsubsection{Constructibility of the intersection of the two focal cones by straightedge (ruler) and compasses}

After the accident Greek we say that if a construction used only a straightedge and compasses, it is called \emph{planar}; if it also required one or more conic sections (other than the circle), then it is called \emph{solid}. The third category included all constructions that did not fall into either of the previous two categories. In this section, we investigate the construction of Chasles from that point of view that the construction of the common edges is planar or not. Of course, the only point which we have to discuss the intersection of those conics from which the first one is the focal hyperbole and the second one is the central projection $\mathcal{E}'$ of the focal ellipse on the plane of the focal hyperbole. The latter conics could be ellipse, parabole, or hyperbole, respectively. In the present subsection we describe the construction of the central projection in that case, when the image is an ellipse. Clearly, if we use analytic description this disjunction is unnecessary.

The equality of the focal ellipse using an appropriate Cartesian coordinate system with axes $x,y,z$ is
$$
\frac{x^2}{a_1^2-a_3^2}+\frac{y^2}{a_2^2-a_3^2}=1 \quad z=0.
$$
The focal hyperbola in the same coordinate system can be given in the form
$$
\frac{x^2}{a_1^2-a_2^2}-\frac{z^2}{a_2^2-a_3^2}=1 \quad y=0.
$$
If the point $O=(x',y',z')^T$ is the centrum of the central projection and $X=(x,y,0)^T$ is a point of the focal ellipse then the parametric equation of the cone through the ellipse is
$$
\gamma(x,y,z)=(x',y',z')^T+t(x-x',y-y',-z')^T.
$$
This cone intersects the plane $y=0$ of the hyperbole in the points $(x'+\frac{y'}{y'-y}(x-x'),0,z'-\frac{y'z'}{y'-y})^T$, where the point $X$ is on the ellipse. (If the plane $y=y'$ intersects the focal ellipse in two points the image has two ideal point (it is a hyperbole) if touches the ellipse then it has one ideal point (it is a parabole) and if do not intersect the ellipse the image has no ideal point (it is an ellipse).) The projection of the point $X$ is on the focal hyperbole if it also satisfies the equality of the hyperbole, hence for the unknown values $x$ and $y$ we have the following system of equations
\begin{eqnarray}\label{eq:commonpoints}
\frac{\left(x'+\frac{y'}{y'-y}(x-x')\right)^2}{a_1^2-a_2^2}-\frac{\left(z'-\frac{y'z'}{y'-y}\right)^2}{a_2^2-a_3^2} & = & 1 \\
\nonumber \frac{x^2}{a_1^2-a_3^2}+\frac{y^2}{a_2^2-a_3^2} & = & 1.
\end{eqnarray}
Assume that $y\ne y'$ and simplify the first equation in (\ref{eq:commonpoints}) ordering it with respect to the decreasing power of $y$. We get that
\begin{equation}\label{eq:orderedfirsteq}
0=\alpha y^2-\beta y- \gamma
\end{equation}
where
$$
\alpha=\left[(a_1^2-a_2^2)(a_2^2-a_3^2)-(a_2^2-a_3^2){x'}^2+(a_1^2-a_3^2){y'}^2+(a_1^2-a_2^2){z'}^2\right]
$$
$$
\beta=\left[2{y'}(a_2^2-a_3^2)\left((a_1^2-a_2^2)-x'x\right)\right] \mbox{ and } \gamma={y'}^2(a_2^2-a_3^2).
$$
We get a planar construction for the required intersection if one of the numbers $x$ and $y$ can be get from the values $a_i^j,x',y',z'$ algebraically using only the addition, subtraction, multiplication, division, and square root operations, since the other one can be get from it by the intersection of a line and an ellipse. In the following situations we get such a constructions:
\begin{itemize}
\item $x'=0$. In this case, $x$ vanish from $\beta$ and $y$ can be constructed from (\ref{eq:orderedfirsteq}). Geometrically this equality means that $OP$ is a profile line in Fig. \ref{fig:imoffocellipse} meaning that the semi diameter $OP$ is orthogonal to the major semi axis of the ellipse determined by the semi diameters $OQ$ and $OR$. In this case, $\alpha=\left[(a_1^2-a_2^2)(a_2^2-a_3^2)+(a_1^2-a_3^2){y'}^2+(a_1^2-a_2^2){z'}^2\right]$ and $\beta^2=4\gamma (a_1^2-a_2^2)^2$ hence we have
    $$
    y_{1,2}=\frac{{y'}(a_2^2-a_3^2)(a_1^2-a_2^2)\pm \sqrt{\gamma \left((a_1^2-a_2^2)^2+\alpha\right)}}{\alpha}=
    $$
    $$
    =\frac{{y'}(a_2^2-a_3^2)(a_1^2-a_2^2)\pm y'\sqrt{(a_2^2-a_3^2)\left((a_1^2-a_2^2)(a_1^2-a_3^2)+ (a_1^2-a_3^2){y'}^2+(a_1^2-a_2^2){z'}^2\right)}}{\left[(a_1^2-a_2^2)(a_2^2-a_3^2)+(a_1^2-a_3^2){y'}^2+(a_1^2-a_2^2){z'}^2\right]}.
    $$
\item $y'=0$. Then the equation is simplified to the equation $0=\alpha y^2$ from which either $\alpha =0$ or $y=0$. Since we assumed that $y'-y\ne 0$ the second case algebraically  cannot be realized here (however it has also geometrical meanings). From $\alpha=(a_1^2-a_2^2)(a_2^2-a_3^2)-(a_2^2-a_3^2){x'}^2+(a_1^2-a_2^2){z'}^2=0$ we can see that $O$ is on the focal hyperbole the searched edges are the lines joint $O$ with the foci of the focal hyperbole and the projection of the focal ellipse to the plane of the focal hyperbole contains only these two lines. The axes of the ellipsoid are the $y$-axis and two lines through the origin parallel to the tangent of the focal hyperbole at $O$ and the normal of the focal hyperbole also at $O$. The lengths of the axis can be determined easily, too. Similar immediate solution can be get for the original construction in the case of any other point of the plane $y=0$ if we choose it for $O$. From the conjugate system of diameters we get this situation when $OP$ is perpendicular to the minor semi axis of the ellipse with conjugate diameters $OR$ and $OQ$.
\end{itemize}

\begin{statement}\label{st:nonconstr}
We can choose the parameters $a:=a_1^2-a_2^2$, $b:=a_2^2-a_3^2$ $x'$, $y'$ and $z'$ on such a way that the common points of the investigated conics can not be constructed by ruler and compasses.
\end{statement}

As a consequence we can see that in this case there is no possibility to give planar construction for this step of the original Chasles's construction, because such a construction will be a construction for the intersection of our conics, too.

\begin{proof}
By the symmetry with respect to the plane $y=0$ we can assume that $y'>0$. First we continue the algebraic determination of the solution of (\ref{eq:orderedfirsteq}). From $a=a_1^2-a_2^2$ and $b=a_2^2-a_3^2$ follows $a+b=a_1^2-a_3^2$. Then we have
$$
\alpha=ab-b{x'}^2+(a+b){y'}^2+a{z'}^2, \quad \beta=2{y'}b(a-{x'}x), \quad \gamma={y'}^2b;
$$
and also
$$
x^2=\sqrt{(a+b)\left(1-\frac{{y}^2}{b}\right)}=\frac{1}{\sqrt{b}}\sqrt{\left(b-{y}^2\right)(a+b)}.
$$
Now we can write
$$
\beta=2{y'}\sqrt{b}\left(a\sqrt{b} \mp {x'}\sqrt{\left(b-{y}^2\right)(a+b)}\right)=\sqrt{\gamma}f(y),
$$
by the function $f(y):=a\sqrt{b} \mp {x'}\sqrt{\left(b-{y}^2\right)(a+b)}$.
A solution of (\ref{eq:orderedfirsteq}) holds one of the equalities
$$
y=\frac{\sqrt{\gamma}}{\alpha}\left(f(y) \pm \sqrt{f(y)^2+\alpha}\right),
$$
equivalently one of the equations
$$
\left(\frac{\alpha}{\sqrt{\gamma}}y-f(y)\right)=\pm \sqrt{f(y)^2+\alpha}.
$$
The square of this equations are
$$
\left(\frac{\alpha}{\gamma}y^2-2\frac{y}{\sqrt{\gamma}}f(y)\right)=1 \quad \mbox{where } f(y)=a\sqrt{b} \mp {x'}\sqrt{\left(b-{y}^2\right)(a+b)}.
$$
From this we get
$$
\left(\frac{\alpha}{\gamma}y^2-\frac{2}{{\sqrt{\gamma}}}a\sqrt{b}y-1\right)^2=\frac{4}{\gamma}y^2{x'}^2b(a+b)-\frac{4}{\gamma}{x'}^2(a+b){y}^4
$$
leading to the equality
$$
\left(\frac{\alpha^2}{\gamma^2}+\frac{4}{\gamma}{x'}^2(a+b)\right){y}^4-4\frac{\alpha a\sqrt{b}}{\gamma\sqrt{\gamma}}y^3+
\left(-2\frac{\alpha}{\gamma}+\frac{4}{\gamma}a^2b - \frac{4}{\gamma}{x'}^2b(a+b)\right)y^2+\frac{4}{{\sqrt{\gamma}}}a\sqrt{b}y+1=0,
$$
where $\alpha=ab-b{x'}^2+(a+b){y'}^2+a{z'}^2$. For $\alpha=0$ we get
$$
\frac{{x'}^2}{{y'}^2}\frac{(a+b)}{b}{y}^4+\frac{\left(a^2 - {x'}^2(a+b)\right)}{{y'}^2}y^2+\frac{a}{y'}y+\frac{1}{4}=0,
$$
and also the assumption
\begin{equation}\label{eq:alpha}
0=ab-b{x'}^2+(a+b){y'}^2+a{z'}^2.
\end{equation}
Using the parameters $a=y'=1$ and $b=x'=2$ we get that $z'=\pm\sqrt{3}$ and the equation simplified to
$$
24{y}^4-44y^2+4y+1=0
$$
Introducing the real numbers $c,d,e$ such that
$$
24{y}^4-44y^2+4y+1=(\sqrt{24}y^2+c)^2-(dy+e)^2=(\sqrt{24}y^2+c+dy+e)(\sqrt{24}y^2+c-dy-e)
$$
holds we can factorize this polynomial into the product of two quadratic terms. Using the method on the comparing of the coefficients, we get for the parameters the following system of equations
$$
2\sqrt{24}c-d^2=-44, \quad de=-2, \quad \mbox{and} \quad 1=c^2-e^2.
$$
Clearly a real root of the original equation can be constructed if and only if the numbers $c,d,e$ can be constructed. This leads to the following cubic equation in $c$
$$
0=2\sqrt{24}c+44-\frac{4}{e^2}=2\sqrt{24}c+44-\frac{4}{c^2-1} \mbox{ equivalently } 0=c^3+\frac{11}{\sqrt{6}}c^2-1c-\frac{12}{\sqrt{6}}.
$$
If $c$ is a root of the latter polynomial than $c +\frac{11}{3\sqrt{6}}$ is a root of the cubic
$$
c^3+\left(-1-\frac{121}{18}\right)c+\left(-\frac{12}{\sqrt{6}}+\frac{11}{3\sqrt{6}}+2\frac{11^3}{6\sqrt{6}27}\right).
$$
Hence we have to investigate the cubic equation
\begin{equation}\label{eq:onc}
0=c^3-\frac{139}{18}c+\frac{328\sqrt{6}}{3^5}.
\end{equation}
The coefficients of it are in $\mathbb{Q}(\sqrt{6})$ which elements are of the form $c=\lambda+\nu\sqrt{6}$ with rational $\lambda $ and $\nu$. Substitute this number into the equality above then we get that it must be hold two equalities $0=\nu^3+3\lambda^2\nu-\frac{139}{18}\nu+\frac{328}{3^5}$ and $0=\lambda^3+18\lambda\nu^2-\frac{139}{18}\lambda=\lambda(\lambda^2+18\nu^2-\frac{139}{18})$. From the second equality we have two possibilities.

First we assume that $\lambda=0$ then the first equality is reduced to the equality $0=\nu^3-\frac{139}{18}\nu+\frac{328}{3^5}$ which has no rational root. In fact, if $\nu=m/n$ where $m$ and $n$ are relative primes, we have the equality
$$
0=2\cdot 3^5m^3+n^2(-3^3\cdot 139 m+2^4\cdot 41 n)=(2\cdot 3^5m^2-3^3\cdot 139 n^2) m+2^4\cdot 41 n^3.
$$
From this it follows that $n=\pm 1, \pm 3, \pm 3^2$. If $n=\pm 1$ then $0=2\cdot 3^5m^3-3^3\cdot 139 m\pm 2^4\cdot 41$ which is impossible because $3 \nmid 2^4\cdot 41$. If $n=\pm 3$ then the equality $0=2\cdot 3^5m^3-3^5\cdot 139 m\pm 2^4\cdot 41\cdot 3^3$ is also wrong because $3^2\nmid 2^4\cdot 41$. In the last subcase the equation $0=2\cdot 3^5m^3-3^7 \cdot 139 m \pm 2^4\cdot 41\cdot 3^5$ is simplified into $0=2m^3-9 \cdot 139 m \pm 2^4\cdot 41$ and thus $m\mid 2^4\cdot 41$ and it is even, hence it is $\pm 2$, $\pm 4$, $ \pm 8$ or $\pm 16$. But these values are not solutions of the equation.

If $\lambda \ne 0$ then $\lambda^2=-18\nu^2+\frac{139}{18}$ and $0=-53\nu^3+ \frac{139}{9}\nu+\frac{328}{3^5}$ which also has no rational roots.
In fact,
$$
0=-53\cdot 3^5m^3+3^3\cdot 139 mn^2+328n^3=-53\cdot 3^5m^3+(3^3\cdot 139 m+328n)n^2
$$
which implies that $n^2=1,3^2$ or $n^2=3^4$. The respective equalities are $0=-53\cdot 3^5m^3+3^3\cdot 139 m\pm 328$, $0=-53\cdot 3^5m^3+3^5\cdot 139 m\pm 328\cdot 3^3$ and $0=-53\cdot 3^5m^3+3^7\cdot 139 m\pm 328 \cdot 3^6$. The first two equations have no integer solutions because $3\nmid 328$ and $3^2\nmid 328$. The last one has the following form $0=-53\cdot m^3+3^2\cdot 139 m\pm 328 \cdot 3$ implying that $m$ has to be divisible by $3$ which is also impossible since $n$ and $m$ have no common divisor.

We shew that the cubic polynomial (\ref{eq:onc}) are irreducible over $\mathbb{Q}(\sqrt{6})$ so it has no constructible roots (see in \cite{szokefalvi}). Since on the other hand it has real roots the statement is fulfill.
\end{proof}

\begin{conjecture}\label{conj:nonconstr}
If $x',y'\ne 0$ we can choose parameters $a:=a_1^2-a_2^2$, $b:=a_2^2-a_3^2$ on such a way that to give planar construction for the determination of the common point of the investigated conics do not possible.
\end{conjecture}

\end{document}